\documentclass[a4paper,11pt]{amsart}
\usepackage[dvips]{graphicx}
\newtheorem{theorem}{Theorem}[section]
\newtheorem{lemma}[theorem]{Lemma}
\newtheorem{proposition}[theorem]{Proposition}
\theoremstyle{definition}
\newtheorem{definition}[theorem]{Definition}
\theoremstyle{remark}
\newtheorem{remark}[theorem]{Remark}
\numberwithin{equation}{section}
\DeclareMathOperator{\Div}{div}
\DeclareMathOperator{\supp}{Supp}
\newcommand{\erre}{\mathbf{R}}
\newcommand{\n}{\textbf{N}}
\newcommand{\dint}{\int\!\!\!\!\int}
\newcommand{\s}{\sigma}
\newcommand{\rs}{\sqrt{\sigma}}
\newcommand{\la}{\lambda}
\newcommand{\be}{\beta}
\newcommand{\frho}{\frac{\rho-\rs}{\rs}}
\newcommand{\fy}{\frac{y}{\rs}}
\newcommand{\ft}{\frac{t-\s}{\s}}
\newcommand{\sud}{\sqrt{x_1^2+x_2^2}}
\newcommand{\z}{\mathbf{z}}
\newcommand{\f}{\mathbf{f}}
\newcommand{\ns}{\mathbf{u}}
\begin{document}
\title{Some examples of singular fluid flows}
\author{Marco Romito}
\address{Dipartimento di Matematica, Universit\`a di Firenze, Viale 
         Morgagni 67/a, 50134 Firenze, Italia}
\email{romito@math.unifi.it}
\subjclass{Primary 76D05; Secondary 35A20}
\keywords{Navier-Stokes equations, emergence of singularities,
          partial regularity, suitable weak solutions}

\begin{abstract}
We explain the construction of some solutions of the Stokes system with
a given set of singular points, in the sense of Caffarelli, Kohn and 
Nirenberg \cite{CKN}. By means of a partial regularity theorem (proved
elsewhere), it turns out that we are able to show the existence of
a suitable weak solution to the Navier-Stokes equations with a 
singular set of positive one dimensional Hausdorff measure.
\end{abstract}
\maketitle
%
%
%

\section{Introduction}

There is a wide interest on the problem of regularity for the Navier-Stokes
equations. One of the most interesting achievement was gained by Caffarelli,
Kohn and Nirenberg in 1982. They showed that a (suitable) weak solution $u$
to the Navier-Stokes equations has a set of singular points of null 
one-dimensional Hausdorff measure. In their definition, a \textsl{regular} 
point for a solution is a space-time point $(t,x)$ such that the solution is
essentially bounded in a small neighbourhood of $(t,x)$. A \textsl{singular}
point is then a space-time point which is not regular. In the sequel, we 
will call \textsl{singular set} the set of all singular points for a weak 
solution $u$.

A few years later Scheffer in \cite{Sc1}, \cite{Sc2} and \cite{Sc3}, found
some examples of solutions to the Navier-Stokes inequality (namely, vector
fields that satisfy only the local energy inequality, but not necessarily
the equations) having a nearly one dimensional singular set. Scheffer pointed
out that the limit of the Caffarelli, Kohn and Nirenberg theorem rested in
the energy method, so that a clever, complete use of the equations could give
a better result.

We are indebted with some of Scheffer's ideas. Since we deal with a linear
equation, our computations are simpler. Nevertheless we ends up with a solution
to the Navier-Stokes equations having a singular set whose dimension is bigger
than the one of the Scheffer's example. This is not in contrast to the
Caffarelli, Kohn and Nirenberg theorem \cite{CKN}, since our body force
is less regular (see Remark \ref{noCKN}).

In view of these considerations, throughout the paper we will call \textsl{thin}
a set having null one-dimensional Hausdorff measure, and \textsl{fat} a set
having positive, possibly infinite, one-dimensional Hausdorff measure. Moreover,
we will denote by $S(u)$ the singular set of a vector field $u$.


\section{Main results}

We will present two different construction of solutions to the 
Stokes equation
\begin{align}\label{stokes}
&\partial_t z-\triangle z+\nabla Q=f\\
&\Div z=0\nonumber
\end{align}
with a given \textsl{fat} singular set. The main idea underlying the
first construction is that the solution has a property of
\textsl{self-similarity}, in a sense related to the one given by Leray 
\cite{Ler}: it is invariant under the scaling
$$
z(t,x)\longrightarrow\la z(\la^2t,\la x).
$$
on a fixed countable set of time intervals. In other words, we start
from a solution to the Stokes equation in a small time
interval. Afterwards, the solution is extended in the next interval by
the self-similarity property, and so on. After a finite time the
singular set of the solution contains a given fractal set.

In the second example we complete the same construction on $\erre^2$,
then we obtain an axially-symmetric solution on $\erre^3$. The main
interest of this example is that it can be used to obtain a solution
to the Navier-Stokes equation with a \textsl{fat} singular set. For this
purpose we use a partial regularity result proved in \cite{FlRo}. The
last part of the section is devoted to the explanation of this result
and of its application on our example.

\subsection{Solutions to the Stokes system with a fat singular set}
The first construction we explain fails to be a solution with a
\textsl{fat} singular set, if at least we want to use it for the
Navier-Stokes equations (see Remark \ref{remfail1}).

Anyway this example has some interest: firstly there is a more
general expression of the body force (it is composed of a mean term
and a fluctuation). Moreover it should be possible to use this
construction in the second example below (see Proposition \ref{trid})
to obtain in that case a even \textsl{fatter} singular set. Finally
the second example exploits some of the idea of the first one, so that
the construction should be plainer. The following proposition is
proved in Section \ref{secflat}.

\begin{proposition}\label{failed}
Let $\alpha\in(0,\frac32)$, then there exist $z_0:\erre^3\to\erre^3$,
$f:\erre_+\times\erre^3\to\erre^3$ and
$g:\erre_+\times\erre^3\to\erre^3$
such that the solution $z$ of the Stokes equation (\ref{stokes})
has the following properties

\begin{itemize}
\item[{\it (i)}] $z\in L^\infty(0,\infty;\left[L^2(\erre^3)\right]^3)\cap L^2(0,\infty;\left[H^1(\erre^3)\right]^3)$,
\item[{\it (ii)}] $z$ has compact support in $[0,\infty)\times\erre^3$,
\item[{\it (iii)}] $z(t)\in C^\infty(\erre^3)$ for each $t\ge 0$,
\item[{\it (iv)}] the singular set $S(z)$ has at least Hausdorff dimension $\alpha$.
\end{itemize}
\end{proposition}

\begin{remark}\label{remfail1}
Unfortunately, the solution given by the proposition above does 
not satisfy the set of assumptions (\ref{assB}) below, so that 
it cannot be used to obtain solutions to the Navier-Stokes
equation with a fat singular set (see Remark \ref{remfail2}).
\end{remark}

From now to the end of this section, we will consider for simplicity
$g\equiv 0$. In order to have a solution of the Stokes system with a fat
singular set satisfying the set of assumptions (\ref{assB}), we need to 
modify the previous construction. The idea is to complete the same 
construction we did for the solution in Proposition \ref{failed} on a 
plane and then to consider an axially-symmetric solution in the space. 
Note that the equation solved by the 2D-solution is different, since 
at the end we want the 3D-solution to solve the Stokes equation. We 
obtain the following result, which will be proved in Section \ref{secaxial}.

\begin{proposition}\label{trid}
There exist $\z_0:\erre^3\to\erre^3$ and $\f:\erre_+\times\erre^3\to\erre^3$ such
that for the solution $\z:\erre_+\times\erre^3\to\erre^3$ of the Stokes system
(\ref{stokes}) the following properties hold
\begin{itemize}
\item[{\it (i)}] $\z\in L^\infty(0,\infty;(L^2(\erre^3))^3)\cap L^2(0,\infty;(H^1(\erre^3))^3)$;
\item[{\it (ii)}] $\z$ has compact support in $[0,\infty)\times\erre^3$;
\item[{\it (iii)}] $\z(t)\in C^\infty(\erre^3)$ for all $t\ge 0$;
\item[{\it (iv)}] $\z$ satisfies the assumptions (\ref{assB}) below;
\item[{\it (v)}] $\mathcal{H}^1(S(\z))>0$.
\end{itemize}
\end{proposition}

\begin{remark}
It should be possible to obtain a \textsl{fatter} singular set than the
one obtained in the previous proposition. The idea is to apply on the
2D solution the same construction as the one given in Proposition 
\ref{failed}, to obtain an axially-symmetric solution with a singular
set of dimension larger than one. Indeed, the example contained
in Proposition \ref{failed} has been given to suggest this
possibility. We have not tried to show this claim because of the
huge amount of annoying computations. Moreover our point is to
show the existence of suitable weak solutions to the Navier-Stokes
equations with a singular set \textsl{fatter} than the one in
Caffarelli, Kohn and Nirenberg \cite{CKN}. This will be given in
the next section.
\end{remark}

\subsection{A suitable weak solution of Navier-Stokes equation with a
fat singular set}

We want to apply now the construction we have given in the previous section 
to the theory of singularities for the Navier-Stokes equations. In \cite{FlRo}
the theory of singularities has been studied in the following way. Consider
the Navier-Stokes equations in a bounded open domain $D\subset\erre^3$, with
regular boundary,
\begin{align}\label{NS}
&\partial_t u-\nu\triangle u+(u\cdot\nabla)u+\nabla P=f+\partial_t g,\notag\\
&\Div u=0,\\
&u(0)=u_0,\notag\\
&u(t,\cdot)=0\qquad\text{on }\partial D,\notag
\end{align}
where the term $f$ represents the mean force and the term $\partial_t g$
represents the fluctuations. The solution of the equation is split
$$
u=v+z,\qquad P=\pi+Q,
$$
where $(z,Q)$ is the solution of the Stokes problem
\begin{align}\label{lNS}
&\partial_t z-\nu\triangle z+\nabla Q=f+\partial_t g,\notag\\
&\Div z=0,\\
&z(0)=0,\notag\\
&z(t,\cdot)=0\qquad\text{on }\partial D,\notag
\end{align}
and $(v,\pi)$ solves the modified Navier-Stokes equations
\begin{align}\label{mNS}
&\partial_t v-\nu\triangle v+\left((v+z)\cdot\nabla\right)(v+z)+\nabla\pi=0,\notag\\
&\Div v=0,\\
&v(0)=u_0,\notag\\
&v(t,\cdot)=0\qquad\text{on }\partial D,\notag
\end{align}

A different notion of suitable weak solution (introduced mostly to treat
the term $\partial_t g$) can be given.

\begin{definition}\label{sws}
Let $(z,Q)$ be the solution of the Stokes problem. A suitable weak 
solution $(u,P)$ to the Navier-Stokes equations (\ref{NS}) is a pair
$$
u\in L^\infty([0,\infty);\left[L^2(D)\right]^3)\cap L^2_{\text{loc}}([0,\infty);\left[H_0^1(D)\right]^3) 
$$
and
$$
P\in L^{5/3}_{\text{loc}}((0,T)\times D)
$$
such that the new variables $v=u-z$ and $\pi=P-Q$ satisfy the modified
Navier-Stokes equations (\ref{mNS}) in the sense of distributions over 
$(0,\infty)\times D$ and moreover satisfy the local energy inequality 
\begin{eqnarray}\label{lei}
\lefteqn{\int_D|v(t)|^2\varphi
         +2\int_0^t\int_D\varphi|\nabla v|^2\le}\nonumber\\
&\le& \int_0^t\int_D|v|^2\left(\frac{\partial\varphi}{\partial t}+\triangle \varphi\right)
      +\int_0^t\int_D\left(|v|^2+2v\cdot z\right)\left((v+z)\cdot\nabla\varphi\right)\\
&   & +2\int_0^t\int_D\varphi z\cdot\left((v+z)\cdot\nabla\right)v
      +\int_0^t\int_D2\pi v\cdot\nabla\varphi\nonumber
\end{eqnarray}
for every smooth function $\varphi:\erre\times D\rightarrow\erre$, 
$\varphi\ge0$, with compact support in $(0,T]\times D$.
\end{definition}

\begin{remark}
For some comments on this definition and on the connection with the
suitable weak solutions, as defined by Caffarelli, Kohn and Nirenberg 
\cite{CKN}, see \cite{FlRo} and \cite{Ro1}.
\end{remark}

The forcing terms $f$, $g$ are taken in the following way
\begin{align}\label{assA}
&f\in L^p_{\text{loc}}((0,\infty)\times D),\qquad p>2,\tag{A}\\
&g\in C^{\frac12-\varepsilon}([0,\infty);H^{2\beta}(D)),\quad g(0)=0,\quad\beta>\varepsilon>0.\notag
\end{align}
The set of assumptions (\ref{assA}) implies (see the appendix in \cite{FlRo})
\begin{align}\label{assB}
&z\in L^\infty_{\text{loc}}([0,\infty);L^2(D))\cap L^2_{\text{loc}}([0,\infty);H^1_0(D)),\notag\\
&z\in L^\infty_{\text{loc}}([0,\infty);L^q_{\text{loc}}(D))\qquad\text{for a }q>6,\tag{B}\\
&\lim_{r\to 0}\frac1r\dint_{Q_r(t,x)}|\nabla z|^2\,ds\,dy=0\qquad\text{for all }(t,x),\notag
\end{align}
where $Q_r(t,x)=\{\,(s,y)\,|\,t-r^2<s<t,\ |x-y|<r\,\}$. Finally, it is shown 
(see Theorem 5.3 in \cite{FlRo}) the following result.

\begin{theorem}
Assume the set of assumptions (\ref{assB}) for $z$. Let $v$ be the solution
of (\ref{mNS}) relative to $z$ and satisfying the local energy inequality
(\ref{lei}). Then
$$
\mathcal{H}^1(S(v))=0
$$
\end{theorem}

In order to give the same conclusion of the previous theorem for the true
solution of Navier-Stokes $u$, we assume the following condition
\begin{align}\label{assC}
&f\in L^p_{\text{loc}}((0,\infty)\times D),\qquad p>\frac52,\tag{C}\\
&g\in C^{\frac12-\varepsilon}([0,\infty);H^{\frac12+2\beta}(D)),\quad g(0)=0,\quad\beta>\varepsilon>0.\notag
\end{align}
In fact the set of assumptions \ref{assC} implies that 
$z\in L^\infty_{\text{loc}}((0,\infty)\times D)$.
\medskip

The idea is now to exploit the gap between the set of assumptions 
(\ref{assA}) and the set of assumptions \ref{assC}, in order to 
find a suitable weak solution of Navier-Stokes equations having
a singular set bigger than it is stated in \cite{CKN}.

Indeed, Proposition \ref{trid} above let us consider a vector field $z$
having a \textsl{fat} singular set and satisfying the set of assumptions
(\ref{assB}). Consequently the solution $v$ of (\ref{mNS}) has a
\textsl{thin} singular set, so that
$$
u=v+z
$$
will have a \textsl{fat} singular set. The complete proof of the following
theorem, which extends the argument above, is given in Section 
\ref{secaxial}

\begin{theorem}\label{wow}
There exist $\ns_0:\erre^3\to\erre^3$, $\f:\erre_+\times\erre^3\to\erre^3$,
a bounded regular domain $D$ and a suitable weak solution $\ns$ of 
the Navier-Stokes system, in the sense of Definition \ref{sws},
such that
$$
\mathcal{H}^1(S(\ns))>0.
$$
\end{theorem}

\begin{remark}
It should be interesting to find which are the minimal assumptions on
the forcing terms (like (\ref{assA}) or (\ref{assC})) which give a
singular set for $z$ of null one-dimensional Hausdorff measure. For
example it may be guessed, with some heuristics, the 
\textsl{critical} summability for $f$. In \cite{Ro2}
it is suggested a way to show this. We know nothing about the term $g$.
\end{remark}


\section{A Stokes vector field with a fractal singular set}\label{secflat}

In this section we prove Proposition \ref{failed}. First we point out
that it is possible to get rid of the divergence-free condition. In 
fact the following construction can be performed separately on the 
three components of a given divergence-free initial condition, giving
a solution $z$ and the body forces $f$, $g$ that will be divergence-free 
vector fields. Consequently there is no loss of generality in
working with a real valued $z$.

For the sake of clarity, the construction is divided into three steps.
In the first step a solution which has a singular point is built, in 
the second step the geometry of the fractal set is explained and 
finally in the third step the arguments of the previous steps are 
glued together to obtain the singular solution.

\subsection{A solution with a singular point}

Consider a function $z_0\in C^\infty_{\text{c}}(\erre^3)$ such that $z_0(0)=1$
and $z_0\ge 0$. Fix two parameters $\s\in(0,1)$ and $\la>1$. Set for $N\ge 1$,
$$
\s_0=0,\quad\s_N=\sum_{k=1}^N\s^k\quad\text{ and }\quad T=\sum_{k=1}^\infty\s^k
$$
and $I_N=[\s_{N-1},\s_N]$. Let $f_1$, $g_1$ be functions in 
$I_1\times\erre^3$ (without loss of generality they can be taken 
smooth and with compact support), with $g_1(0)=0$,  such that the 
solution $z_1$ of the Stokes equation is positive, smooth, 
with compact support and
\begin{eqnarray*}
&z_1(0,x)=z_0(x),\qquad z_1(\s,x)=\la z_0(\frac{x}{\sqrt\s})\\
&\|z_1(t)\|_{L^p(\erre^3)}\le C\|z_0\|_{L^p(\erre^3)}\qquad\text{for each }p\ge 2\text{ and }t\ge0. 
\end{eqnarray*}
Then for $N\ge 2$, functions $z_N$, $f_N$ and $g_N$ are defined 
inductively in $I_N\times\erre^3$ as
\begin{align*}
&z_N(t,x)=\la z_{N-1}(\frac{t-\s}{\s},\frac{x}{\sqrt{\s}}),\\
&f_N(t,x)=\frac\la\s f_{N-1}(\frac{t-\s}{\s},\frac{x}{\sqrt{\s}}),\\
&g_N(t,x)=g_{N-1}(\s_{N-1})+\la\left[g_{N-1}(\frac{t-\s}{\s},\frac{x}{\rs})-g_{N-1}(\s_{N-2},\frac{x}{\rs})\right].
\end{align*}
Finally set
$$
z(t,x)=z_N(t,x),\qquad f(t,x)=f_N(t,x)\quad\text{and}\quad g(t,x)=g_N(t,x)
$$
if $t\in I_N$ and $z(t,x)=f(t,x)=g(t,x)=0$ if $t\ge T$.

It is easy to verify that, if
$$
\la\s^{3/4}<1,
$$
then $z\in L^2(0,\infty;L^2(\erre^3))\cap L^2(0,\infty;H^1(\erre^3))$ 
and $z$ solves the equation in the sense of distributions. Finally it can 
be verified that
$$
z(\s_N,0)=z_N(\s_N,0)=\la^Nz_0(0)=\la^N
$$
and so the point $(T,0)$ is a singular point for $z$.

\begin{remark}
It is easy to see that $f\in L^p(\erre_+\times\erre^3)$ 
if
$$
\la^p\s^{\frac52-p}<1,
$$
while $g\in C^{\frac12-\varepsilon}([0,\infty);H^{2\beta}(D))$ if
$$
\la\s^{\frac14-(\beta-\varepsilon)}\le 1.
$$
\end{remark}

\subsection{The geometrical construction}

A classical construction of self similar fractal sets will be used 
and it will give a generalised Cantor set (see for example David,
Semmes \cite{DaSe}). Fix two integers $k$ and $m\le k^3$ and 
consider the unit cube $Q$ in $\erre^3$. The first step of the 
construction is to divide the cube in $k^3$ equal cubes and to 
choose $m$ of them.

The second step is to perform the same operation on each of these 
$m$ cubes: take each of them, divide it in $k^3$ parts and select 
$m$ of these smaller cubes in the same relative positions as in 
the first step. Now there are $m^2$ cubes, the next step is to 
apply this transformation to each of them to obtain $m^3$ cubes 
and so on (see the figure which shows the first three iterations 
with $k=5$ and $m=8$).

\begin{figure}[ht]
\includegraphics[width=3.91cm]{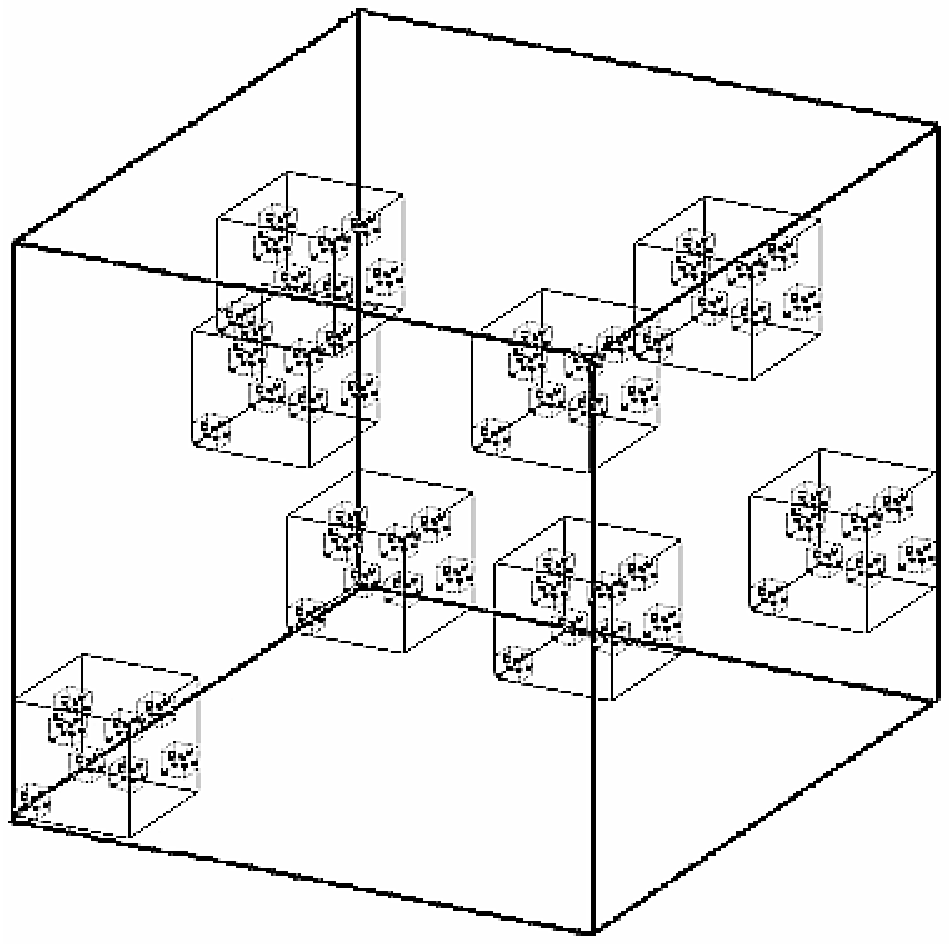}
\end{figure}

The fractal set $C_{k,m}$ will be the only compact subset of 
$\erre^3$ which is invariant under the following transformation. 
If $C_{k,m}$ is re-scaled with a factor $k^{-1}$, $m$ copies of 
the re-scaled set are made and each copy is put in the place of 
each of the cubes chosen, then $C_{k,m}$ is obtained again. So 
it is easy to calculate the Hausdorff dimension $\alpha$ of 
$C_{k,m}$, in fact
$$
\frac{m}{k^{\alpha}}=1,\quad\text{that gives }\alpha=\frac{\log m}{\log k}
$$
(see for example David and Semmes \cite{DaSe}).

\subsection{A fractal set of singularities}

The construction of the fractal set and the construction of the 
first step are mixed up together: there the starting point was 
a solution in the first interval $I_1$ that was $z_0$ in the 
beginning and a re-scaling of $z_0$ at the end. Now the basic 
solution will be constructed as a function which solves the
equation, it is equal to $z_0$ at time $t=0$ and is equal to 
$m$ re-scaling of $z_0$ placed in the points corresponding to 
point $x=0$ in the construction of the fractal set. At each 
step the construction will be iterated, following the outline 
of the procedure of the construction of the fractal set. 

More precisely, fix integers $k$ and $m$ as above and set 
$\s=k^{-2}$ and $\s_N$ and $T$ as in the previous paragraph. 
Consider the $m$ points $x_1$, $x_2$, \dots, $x_m$ vertices 
of the $m$ chosen small cubes which are in the same position 
as $x=0$ in the bigger cube and set
$$
\be_i(x)=k(x-x_i)=\frac{x-x_i}{\sqrt{\s}}\qquad i=1,\ldots,m
$$
and for each function $u=u(t,x)$,
$$
\left(S_i(u)\right)(t,x)=u(\frac{t-\sigma}{\sigma},\be_i(x)).
$$
Fix $i\in\{1,\ldots,m\}$ and define functions $Z_i$, $F_i$
and $G_i$, with $G_i(0)=0$, in such a way that $Z_i$ is the 
solution of the equation with forcing term $F_i+\partial_t G_i$ and
$$
Z_i(0,x)=\frac1m z_0(x)\qquad Z_i(\s,x)=\frac\la{m}z_0(\be_i(x)),
$$
with $Z_i$ enjoying the same properties as the basic solution 
of the first part.

For each integer $N\ge 2$ and for all 
$i\in\{1,\ldots,m\}$ we define in $I_N\times\erre^3$ the
following functions
$$
Z_N=\sum_{i=1}^m(Z_{N-1})_i,
\qquad
F_N=\sum_{i=1}^m(F_{N-1})_i
\quad\text{and}\quad
G_N=\sum_{i=1}^m(G_{N-1})_i,
$$
where
\begin{eqnarray*}
(Z_{N-1})_i&=&\frac\la{m}S_i(Z_{N-1}),\\
(F_{N-1})_i&=&\frac\la{m\s}S_i(F_{N-1}),\\
(G_{N-1})_i(t)&=&\frac1mG_{N-1}(\s_{N-1})+\frac\la{m}\left[S_i(G_{N-1})(t)-S_i(G_{N-1})(\s_{N-1})\right].
\end{eqnarray*}
Finally define
$$
z(t,x)=Z_N(t,x)\qquad\text{if }t\in I_N
$$
and $z(t,x)=0$ if $t\ge T$. In a similar way $f$ and $g$ are defined.

With some calculations (that we do not complete, since are similar, but
easier, to the calculations of the next section), the following proposition 
can be obtained.

\begin{proposition}
The following properties hold
\begin{itemize}
\item[(i)]  for each $N\ge1$, $z_N(\s_{N-1})=z_{N-1}(\s_{N-1})$ and
$$
z_{N-1}(\s_{N-1},x)=\sum_{i_1,\ldots,i_{N-1}=1}^m\left(\frac\la{m}\right)^{N-1}z_0(\be_{i_1}\circ\ldots\be_{i_N}(x)),
$$
\item[(ii)] if $\la\s< 1$, $z$ is a solution in the sense of distributions of 
the equation with body force $f+\partial_tg$,
\item[(iii)] if $\la\s^{\frac34}<1$ then $z\in L^\infty(0,\infty;L^2(\erre^3))\cap L^2(0,\infty;H^1(\erre^3))$.
\end{itemize}
\end{proposition}

Finally, the fractal set constructed in the previous paragraph is contained
in the singular set of $z$. Indeed, define subsets $A_0=\{0\}\subset\erre^3$, 
$A_1=\{x_1,\ldots,x_N\}$ and for arbitrary $N$
$$
A_N=\{x\in\erre^3\,|\,\text{there exist }i_1,\ldots,i_N\text{ such that }\be_{i_1}\circ\ldots\circ\be_{i_N}(x)=0\,\}.
$$
Each point of the fractal set $C_{k,m}$ is the limit of a sequence 
$(y_N)_{N\in\n}$, with $y_N\in A_N$, so it is sufficient to show that
$z(\s_N,x)$, with $x\in A_N$, goes to infinity as $N$ goes to infinity. 
In fact let $x\in A_N$ and let $j_1,\ldots,j_N$ be such that 
$\be_{j_1}\circ\ldots\circ\be_{j_N}(x)=0$, then
\begin{eqnarray*}
z(\s_N,x)=z_N(\s_N,x)&=&\sum_{i_1,\ldots,i_N=1}^m\left(\frac\la{m}\right)^Nz_0(\be_{i_1}\circ\ldots\circ\be_{i_N}(x))\\
&\ge&\left(\frac\la{m}\right)^Nz_0(\be_{j_1}\circ\ldots\circ\be_{j_N}(x))=\left(\frac\la{m}\right)^N,
\end{eqnarray*}
since $z_0$ is positive. In conclusion if $\la>m$, the claim is true.

The condition $\la\s^{3/4}<1$ implies that $z$ is a solution of finite 
energy, and this forces the Hausdorff dimension of the singular set 
to be any number less than $\frac32$ since
$$
{\dim}_{\mathcal H}C_{k,m}=\frac{\log m}{\log k}<\frac{\log\la}{\log k}=2\frac{\log\la}{-\log\s}<2\cdot\frac34=\frac32,
$$
and $\la$ can be arbitrarily close to $m$.

\begin{remark}\label{remfail2}
Again, we can check that if $\la^p\s^{\frac52-p}<1$ then 
$f\in L^p((0,\infty)\times\erre^3)$ and, if 
$\la\s^{\frac14-(\beta-\varepsilon)}\le 1$ then
$g\in C^{\frac12-\varepsilon}([0,\infty);H^{2\beta}(D))$.

Moreover $z\in L^\infty(0,\infty;L^q(D))$ if $\la^q\s^{\frac32}<1$.
Therefore, in view of assumptions (\ref{assB}), we can give an 
estimate on the Hausdorff dimension of the set $C_{k,m}$ in 
dependence of $q$:
$$
\dim_{\mathcal{H}}C_{k,m}<\frac3q,
$$
so that if $q>6$, the dimension is smaller than $\frac12$.
\end{remark}


\section{The axially-symmetric construction}\label{secaxial}

In this section we will prove Proposition \ref{trid} and Theorem
\ref{wow}. The first lemma gives the axially-symmetric construction
for the solution of equation (\ref{heat}). The Proposition
\ref{trid} will be proved once we use the lemma to obtain a
divergence-free vector field enjoying the same properties. In the 
second part of the section, Theorem \ref{wow} is proved. 

\begin{lemma}\label{oned}
There exist $Z_0:\erre^3\to\erre$ and $F:\erre_+\times\erre^3\to\erre$ such
that for the solution $Z:\erre_+\times\erre^3\to\erre$ of the equation
\begin{align}\label{heat}
&\partial_t Z -\triangle Z=F\qquad\text{in }\erre_+\times\erre^3,\\
&Z(0)=Z_0\notag
\end{align}
the following properties hold
\begin{itemize}
\item[{\it (i)}] $Z\in L^\infty(0,\infty;L^2(\erre^3))\cap L^2(0,\infty;H^1(\erre^3))$;
\item[{\it (ii)}] $Z$ has compact support in $[0,\infty)\times\erre^3$;
\item[{\it (iii)}] $Z(t)\in C^\infty(\erre^3)$ for all $t\ge 0$;
\item[{\it (iv)}] $Z$ satisfies the assumptions (\ref{assB});
\item[{\it (v)}] $\mathcal{H}^1(S(Z))>0$.
\end{itemize}
\end{lemma}

\begin{proof}
We consider $\la$, $\s$, $\s_N$ and $I_N$ as in the previous section,
we assume
$$
\la^q\s<1\qquad\text{for a }q>6
$$
and we suppose also that $\s<\frac14$. Moreover we set
$$
\rho_N=\sum_{k=1}^N(\rs)^k\qquad\text{and}\qquad\rho_\infty=\sum_{N=1}^\infty(\rs)^N.
$$
Let $z_0=z_0(\rho,y)\in C^\infty_c(\erre^2)$ be a function such that
\begin{eqnarray}\label{zetazero}
\supp z_0\subset(-1,1)\times(-M,M) &\qquad& z_0(0,0)=1\\
\frac{\partial^k z_0}{\partial\rho^k}(t,0,y)=0,\quad\text{for }k=1,2,\nonumber
\end{eqnarray}
for some $M>0$. Then there exist $z_1:\erre_+\times\erre^2\to\erre$
and $f_1:\erre_+\times\erre^2\to\erre$ such that
\begin{equation}\label{ctre}
\frac{\partial z_1}{\partial t}-\triangle z_1-\frac1\rho\frac{\partial z_1}{\partial\rho}=f_1
\qquad t\in I_1,\quad (\rho,y)\in\erre^2.
\end{equation}
Moreover
$$
z_1(0,\rho,y)=z_0(\rho,y)\qquad\text{and}\qquad z_1(\s,\rho,y)=\la z_0(\frho,\fy)
$$
and
$$
\int_{\erre^2}|z_1(t)|^q\le C(q)\qquad\text{for all }t\in I_1\quad\text{and each}\quad q\ge 1.
$$
If $N\ge 2$ and $t\in I_N$, we set for $(\rho,y)\in\erre^2$,
\begin{equation}\label{cuno}
z_N(t,\rho,y)=\lambda z_{N-1}(\ft,\frho,\fy)
\end{equation}
and
$$
f_N(t,\rho,y)=\frac\la\s f_{N-1}(\ft,\frho,\fy)
$$
Finally we set
$$
Z(t,x)=\begin{cases}
       z_N(t,\sud,x_3)&\text{if } t\in I_N,\\
       0              &\text{if } t\ge T,
       \end{cases}
$$
and
$$
F(t,x)=\begin{cases}
       f_N(t,\sud,x_3)&\text{if } t\in I_N,\\
       0              &\text{if } t\ge T,
       \end{cases}
$$

First, we observe that $(\s_N,\rho_N,0)\to(T,\rho_\infty,0)$ and
$$
z_N(\s_N,\rho_N,0)=\ldots=\la^Nz_0(0,0)=\la^{N},
$$
so that
$$
\{\,(t,x)\in\erre_+\times\erre^3\,|\,t=T,\ x_1^2+x_2^2=\rho_\infty^2,\ x_3=0\,\}\subset S(Z)
$$
and property \textit{(v)} holds true.

Then we prove property \textit{(ii)}. By assumptions (\ref{zetazero}) we
can deduce that $\supp z_1(t)\subset(-1,1)\times(-M,M)$, so that,
using (\ref{cuno}), for $t\in I_N$,
\begin{equation}\label{cdue}
\supp z_N(t)\subset (\rho_{N-1},\rho_N+\s^{N/2})\times(-\s^{\frac{N-1}{2}}M,\s^{\frac{N-1}{2}}M),
\end{equation}
and property \textit{(ii)} holds true.

Property \textit{(iii)} is obvious. We prove property \textit{(i)}. First 
we show that $Z\in L^2(0,\infty;H^1(\erre^3))$. Note that, if $t\in I_N$,
$$
|\nabla Z(t,x)|=|\nabla z_N(t,\sud,x_3)|
$$
and, using (\ref{zetazero}), this is true for all $x\in\erre^3$. Then by
the change of variable $x_1=\rho\cos\theta$, $x_2=\rho\sin\theta$ and $x_3=y$,
$$
\int_0^{+\infty}\int_{\erre^3}|\nabla Z(t,x)|^2\,dx\,dt
= 2\pi\sum_{N=1}^\infty\int_{I_N}\int_\erre\int_0^\infty\rho|\nabla z_N(t,\rho,y)|^2\,d\rho\,dy\,dt
$$
and, by formula (\ref{cuno}) and iterating a change of variables $N-1$ times,
\begin{eqnarray*}
\lefteqn{\int_{I_N}\int_\erre\int_0^\infty\rho|\nabla z_N(t,\rho,y)|^2\,d\rho\,dy\,dt=}\\
\qquad&=& \la^2\s \int_{I_{N-1}}\int_\erre\int_{-1}^\infty(\rho_1+\rs\rho)|\nabla z_{N-1}(t,\rho,y)|^2\,d\rho\,dy\,dt=\ldots=\\
&=&       (\la^2\s)^{N-1} \int_{I_1}\int_\erre\int_{-\frac{\rho_{N-1}}{\s^{\frac{N-1}2}}}^\infty(\rho_{N-1}+\s^{\frac{N-1}2}\rho)|\nabla z_1(t,\rho,y)|^2\,d\rho\,dy\,dt\\
&\le&     (\la^2\s)^{N-1} (1+\rho_\infty)\int_{I_1}\int_{\erre^2}|\nabla z_1(t,\rho,y)|^2\,d\rho\,dy\,dt,
\end{eqnarray*}
where the last bound comes from formula (\ref{cdue}).

In the same way, it is possible to show that $Z\in L^\infty(0,\infty;L^q(\erre^3))$.

In order to show property \textit{(iv)}, it is sufficient to show that
$$
\lim_{N\to\infty}\frac1{r_N}\dint_{Q_{r_N}}|\nabla Z|^2=0
$$
for $r_N=\sqrt{T\s^N}$ and a cylinder $Q_{r_N}$ centred in a 
point $(T,x^0)$ (for the other points the property is obvious). In 
fact, by similar computations as in the proof of property \textit{(i)}, 
we obtain
$$
\int_{\s_N}^T\int_{B_{r_N}}|\nabla Z|^2\,dx\,dt
=2\pi\sum_{k=N+1}^\infty\int_{I_{N+1}}\int_{x_3^0-r_N}^{x_3^0+r_N}\int_0^{r_N}\rho|\nabla z_{N+1}|^2\,d\rho\,dy\,dt
$$
and so, by $N$ changes of variables and using (\ref{cdue}),
$$
\int_{I_{N+1}}\int_{x_3^0-r_N}^{x_3^0+r_N}\int_0^{r_N}\rho|\nabla z_{N+1}|^2\,d\rho\,dy\,dt
\le  C(\la^2\s)^N,
$$
where $C$ is a constant independent of $N$. Finally, if 
$r\in(\sqrt{T\s^{N+1}},\sqrt{T\s^N})$,
$$
\frac1r\int_{T-r^2}^T\int_{B_r(x^0)}|\nabla Z|^2
\le \frac1{r_{N+1}}\int^T_{T-r_N^2}\int_{B_{r_N}(x^0)}|\nabla Z|^2\le C(\la^2\rs)^N\to 0.
$$

In order to conclude the proof of the proposition, we have only to show
that $Z$ is a solution in the sense of distributions of equation (\ref{heat}). Let 
$\phi\in C^\infty_c(\erre_+\times\erre^3)$, we have to show that
$$
\int_0^\infty\int_{\erre^3}(Z\frac{\partial\phi}{\partial t}+Z\triangle\phi+\phi F)\,dx\,dt
+\int_{\erre^3}Z(0)\phi(0)\,dx=0.
$$
Since $Z$ solves equation (\ref{ctre}), it follows that
\begin{eqnarray*}
\int_0^\infty\int_{\erre^3}(Z\frac{\partial\phi}{\partial t}+Z\triangle\phi+\phi F)\,dx\,dt
&=& \sum_{N=1}^\infty\left[\int_{\erre^3}Z(t)\phi(t)\,dx\right]^{t=\s_N}_{t=\s_{N-1}}\\
&=& \lim_{N\to\infty}\int_{\erre^3}Z(\s_N)\phi(\s_N)\,dx.
\end{eqnarray*}
This limit is equal to zero since, by the usual \textsl{formula (\ref{cuno})
plus $N-1$ changes of variables plus formula (\ref{cdue})} argument, 
it follows that
$$
\int_{\erre^3}Z(\s_N)\phi(\s_N)\,dx
\le\ldots\le
C \|\phi\|_{L^\infty}(\la\s)^N,
$$
where the constant $C$ does not depend on $N$.
\end{proof}

\subsection{The proof of Proposition \ref{trid}} So far, we have 
considered real valued solutions of the equation. In order to pass to
divergence-free vector fields, we consider the functions $Z$, 
$F$ obtained in the previous lemma. We ask also that
\begin{equation}\label{cquattro}
\int_\erre Z(0,x_1,x_2,\xi)\,d\xi=0\qquad\text{for all }(x_1,x_2)\in\erre^2
\end{equation}
(this reduces to a similar assumption on the function $z_0$ of the previous
lemma). We define
$$
\z(t,x)=\left(Z(t,x),Z(t,x),\z_3(t,x)\right)
\quad\text{and}\quad
\f(t,x)=\left(F(t,x),F(t,x),\f_3(t,x)\right)
$$
where $\z_3$ is defined as
$$
\z_3(t,x)=-\int_{-\infty}^{x_3}\left[\frac{\partial Z}{\partial x_1}(t,x_1,x_2,\xi)+\frac{\partial Z}{\partial x_2}(t,x_1,x_2,\xi)\right]\,d\xi,
$$
in such a way that $\z$ is divergence free. The function $\f_3$ will be 
defined later. Consequently the proposition is proved if we show that 
properties \textit{(i)-(iv)} of Lemma \ref{oned} hold for $\z_3$.
Without loss of generality, we consider
$$
\z_3(t,x)=-\int_{-\infty}^{x_3}\frac{\partial Z}{\partial x_1}(t,x_1,x_2,\xi)\,d\xi.
$$

Property \textit{(ii)} follows from assumption (\ref{cquattro}), while 
\textit{(iii)} is obvious. We prove property \textit{(i)}. We start by
proving that $\z_3\in L^\infty(0,\infty;H^1(\erre^3))$. By direct
computation
$$
\frac{\partial\z_3}{\partial x_1}=-\int_{-\infty}^{x_3}\frac{\partial^2Z}{{\partial x_1}^2}\,d\xi,
\qquad
\frac{\partial\z_3}{\partial x_2}=-\int_{-\infty}^{x_3}\frac{\partial^2Z}{\partial x_1\partial x_2}\,d\xi,
\qquad
\frac{\partial\z_3}{\partial x_3}=-\frac{\partial Z}{\partial x_1}.
$$
Obviously the third term is in $L^2$, so we consider only the first
derivative (the second can be treated in the same way). Again a
direct computation gives
\begin{eqnarray*}
\frac{\partial\z_3}{\partial x_1}
&=& -\int_{-\infty}^{x_3}\frac{x_2^2}{(x_1^2+x_2^2)^{3/2}}
                         \frac{\partial z_N}{\partial\rho}(t,\sud,\xi)\,d\xi\\
& & -\int_{-\infty}^{x_3}\frac{x_1^2}{x_1^2+x_2^2}
                         \frac{\partial^2 z_N}{{\partial\rho}^2}(t,\sud,\xi)\,d\xi,
\end{eqnarray*}
for $t\in I_N$. Note that this formula holds true also when $x_1=x_2=0$,
by virtue of (\ref{zetazero}). For the sake of simplicity, we examine
only the second term. We proceed as in the proof of Lemma
\ref{oned}: we divide the integral in time in an infinite sum and we
estimate each term.
\begin{eqnarray*}
\lefteqn{\int_{I_N}\int_{\erre^3}\left|\int_{-\infty}^{x_3}\frac{x_1^2}{x_1^2+x_2^2}
                         \frac{\partial^2 z_N}{{\partial\rho}^2}(t,\sud,\xi)\,d\xi,\right|^2\,dx}\\
&=& \int_{I_N}\int_\erre\int_0^\infty\int_0^{2\pi}\rho\left|\int_{-\infty}^{y}\cos^2\theta
                         \frac{\partial^2 z_N}{{\partial\rho}^2}(t,\rho,\xi)\,d\xi,\right|^2=\ldots=\\
&=& C\la^2\s\int_{I_{N-1}}\int_\erre\int_{-1}^{+\infty}(\rs\rho+\rs)\left|\int_{-\infty}^y
                         \frac{\partial^2 z_{N-1}}{{\partial\rho}^2}(t,\rho,\xi)\,d\xi,\right|^2=\ldots=\\
&=& C(\la^2\s)^{N-1}\int_{I_1}\int_\erre\int_{-\frac{\rho_{N-1}}{\s^{\frac{N-1}2}}}^{+\infty}(\s^{\frac{N-1}2}\rho+\rho_{N-1})\left|\int_{-\infty}^y
                         \frac{\partial^2 z_1}{{\partial\rho}^2}(t,\rho,\xi)\,d\xi,\right|^2\\
&\le& C (\la^2\s)^{N-1},
\end{eqnarray*}
where $C$ does not depend on $N$. In the same way it is possible to prove
that $\z_3\in L^\infty(0,\infty;L^q(\erre^3))$.

Moreover a modification of the previous computation can be used to show that
also property \textit{(iv)} holds: it is sufficient to proceed as in the
proof of the previous lemma.

Finally we show that $\z_3$ solves equation (\ref{heat}). We set
$$
\f_3(t,x)=-\int_{-\infty}^{x_3}\left(\frac{\partial F}{\partial x_1}+\frac{\partial F}{\partial x_2}\right)(t,x_1,x_2,\xi)\,d\xi.
$$
It is easy to see then that in each interval $I_N$, $\z_3$ solves the equation
$$
\frac{\partial \z_3}{\partial t}-\triangle\z_3=\f_3.
$$
Consequently it is sufficient, as in the proof of the previous lemma, to
show that for each test function $\phi$,
$$
\lim_{N\to\infty}\int_{\erre^3}\z_3(\s_N)\phi(\s_N)=0,
$$
and this can be done in the same way as before. The proof of the
proposition is completed.

\begin{remark}\label{noCKN}
It is easy to see that $\f\not\in L^2((0,\infty)\times\erre^3)$. In fact,
$$
\int_0^\infty\int_{\erre^3}|\f_1(t,x)|^p\,dx\,dt
=\sum_{N=1}^\infty\int_{I_N}\int_{\erre^3}|\f_1(t,x)|^p\,dx\,dt
$$
and
$$
\int_{I_N}\int_{\erre^3}|\f_1(t,x)|^p\,dx\,dt=\ldots=C(\la^p\s^{2-p})^N,
$$
so that $p<2$. More precisely
$$
p\le\frac{2q}{1+q}.
$$
\end{remark}

\subsection{The proof of Theorem \ref{wow}} Now we can conclude the 
proof of Theorem \ref{wow}. From the previous proposition we have a 
vector field $\z$ solution of the Stokes system, with initial condition
$\z(0)\in C^\infty_{\text{c}}$. Let $\overline\z$ be the solution of 
the Stokes system with initial condition $-\z(0)$ and a suitable regular
body force in such a way that $\overline\z$ has compact support in space.
We set $\zeta=\z+\overline\z$. The vector field $\overline\z$ is regular,
so that
$$
S(\zeta)=S(\z).
$$
Finally, $\zeta$ is a solution of the Stokes system with zero initial
condition, so that if $v$ is the solution of the modified Navier-Stokes
equation (\ref{mNS}) associated to $\zeta$ in a suitable regular bounded
domain which contains the support of $\zeta$ (the existence of such solution
is stated in \cite{Ro1}), the vector field
$$
\ns=v+\zeta
$$
is a suitable weak solution of Navier-Stokes equation and
$$
\mathcal{H}^1(S(\ns))>0.
$$


\bibliographystyle{amsplain}

\end{document}